\documentclass[]{interact}

\usepackage{epstopdf}% To incorporate .eps illustrations using PDFLaTeX, etc.
\usepackage{caption}
\usepackage{subcaption}

\usepackage[numbers,sort&compress]{natbib}% Citation support using natbib.sty
\bibpunct[, ]{[}{]}{,}{n}{,}{,}% Citation support using natbib.sty

\theoremstyle{plain}% Theorem-like structures provided by amsthm.sty
\newtheorem{theorem}{Theorem}[section]
\newtheorem{lemma}[theorem]{Lemma}

\newtheorem{proposition}[theorem]{Proposition}

\theoremstyle{definition}

\theoremstyle{remark}

\DeclareMathOperator{\sign}{sign}

\begin{document}

{\centering%
\noindent\begin{tabular}{|l|}
%\centering
\hline
Submitted to Mathematical and Computer Modelling of Dynamical Systems.\\
\hline
\end{tabular}
}%

\articletype{ARTICLE TEMPLATE}% Specify the article type or omit as appropriate

\title{Discrete-time staged progression epidemic models}

\author{
\name{Luis Sanz-Lorenzo\textsuperscript{a}\thanks{CONTACT Author. Email: luis.sanz@upm.es} and Rafael Bravo de la Parra \textsuperscript{b}}
\affil{\textsuperscript{a}Depto. Matem\'aticas, E.T.S.I Industriales, Technical University of Madrid, Madrid, Spain; \textsuperscript{b}U.D. Matem\'aticas, Universidad de Alcal\'a, Alcal\'a de Henares, Spain.}
}

\maketitle

\begin{abstract}
In the Staged Progression (SP) epidemic models, infected individuals are classified into a suitable number of states. The goal of these models is to describe as closely as possible the effect of differences in infectiousness exhibited by individuals going through the different stages.
The main objective of this work is to study, from the methodological point of view, the behavior of solutions of the discrete time SP models without reinfection and with a general incidence function. Besides calculating $\mathcal{R}_{0}$, we find bounds for the epidemic final size, characterize the asymptotic behavior of the infected classes, give results about the final monotonicity of the infected classes, and obtain results regarding the initial dynamics of the prevalence of the disease. Moreover, we incorporate into the model the probability distribution of the number of contacts in order to make the model amenable to study its effect in the dynamics of the disease.
\end{abstract}


\section{Introduction}

\label{sec1}

The appearance and spread of COVID-19 has led to the proposal and analysis of
numerous mathematical models. The simplest are in the form of compartmental
models that distinguish a few states of individuals with respect to disease.
Classic representatives of these models are the SI, SIR and SEIR models, along
with their versions supporting reinfection \cite{brauer_mathematical_2019}. A
natural extension of these models are the so-called Staged Progression (SP)
models \cite{hyman_differential_1999}. In them, infected individuals are
classified into a suitable number of states through which they progress. The
goal of these models is to describe as closely as possible the effect of
differences in infectiousness exhibited by individuals going through the
different stages. A paradigmatic case frequently analyzed through SP models is
that of HIV transmission, see references in \cite{guo_global_2006}.

The main objective of this work is to study, from the methodological point of
view, the behavior of the solutions of the discrete time SP models without reinfection.

In the last years discrete-time epidemic models are receiving increasing
attention
\cite{brauer_discrete_2010,bravo_de_la_parra_reduction_2023,hernandez-ceron_discrete_2013,
hernandez-ceron_reproduction_2013,kreck_back_2022,van_den_driessche_demographic_2019,
van_den_driessche_disease_2019,yakubu_discrete-time_2021}. In
\cite{diekmann_discrete-time_2021} a discrete-time version of the
Kermack-McKendrick continuous-time epidemic model of 1927
\cite{kermack_contribution_1927} is presented, and discrete-time models are
motivated. One of the advantages of discrete-time models is their adaptation
to census data, which are usually collected in regular periods of time. The
other great advantage is their direct numerical implementation, which greatly
facilitates all kinds of simulations.

Most of the works that treat SP models in continuous time do so including
a basic demographic turnover. Apart from the already mentioned \cite{hyman_differential_1999} and \cite{guo_global_2006}, in \cite{guo_global_2008} it is presented an SP model that allows for the amelioration of infected individuals, in \cite{hyman_epidemic_2009} models with differential susceptibilities and staged-progressions are analyzed, and
a general class of multistage epidemiological models that allow for the possible deterioration and amelioration between any two infected stages are studied in \cite{guo_global_2012}.  There is also a work in which the SP model presented, in addition to demographic turnover, includes reinfection, \cite{melesse_global_2010}. These works must be considered models for endemic diseases. Their analyses focus on the existence and stability of equilibria. The basic reproduction number $\mathcal{R}_{0}$ determines, when it is less than 1, that
the disease-free equilibrium (DFE) is globally asymptotically stable, i.e.,
the disease is extinguished. Otherwise, when $\mathcal{R}_{0}>1$, there is a
positive equilibrium, the endemic equilibrium (EE), which is globally
asymptotically stable, which ensures the endemicity of the disease.

In the epidemic models we deal with in this work the interval of study is
short enough that demographic effects need not be taken into account, and
infected individuals are assumed to acquire complete immunity
\cite{brauer_mathematical_2017}. Loosely speaking, the role of $\mathcal{R}%
_{0}$ in these models is to distinguish between the infection dying out and
the onset of an epidemic. The final size of the epidemic is the main
characteristic to highlight. In continuous time we can mention two references
on SP epidemic models, \cite{brauer_discrete_2010,brauer_simple_2011}. In
\cite{brauer_discrete_2010}, a general SP epidemic model is built and to
analyze its behavior a relation between the initial and final sizes of the
number of susceptibles and $\mathcal{R}_{0}$ is found. The above is extended
to treatment models. In \cite{brauer_simple_2011}, the effect of individual
behavioral changes, essentially reducing the number of contacts, on the final
size of the epidemic is analyzed. This is first done on a SIR model and then
generalized to an SP model.

The aforementioned \cite{brauer_discrete_2010} is the seminal work on
discrete-time SP epidemic models. These models are presented for a particular
type of incidence function, the value of $\mathcal{R}_{0}$ is obtained and, as
in the continuous case, the equation that relates $\mathcal{R}_{0}$ and the
final size of the epidemic is found. In this paper we deal with discrete-time
SP epidemic models with a general incidence and develop a detailed analysis of
the behaviour of solutions. In particular, we find bounds for the final size
of the epidemic, characterize the speed of convergence of the infected classes
to zero, give results about the final monotonicity of the infected classes,
and obtain results regarding the behavior of the prevalence of the disease,
with special attention to conditions under which the prevalence initially
rises or decays monotonically to zero.

The structure of the paper is as follows: In Section \ref{sec2}, a
discrete-time staged progression epidemic model with general incidence, no
demography, and no re-infection, is presented and some of its basic properties
are analyzed, including the $\mathcal{R}_{0}$ calculation. In Section
\ref{sec3} several behavioral results of the long-term behaviour of solutions
of the model are presented. Among others, the existence of a limit
distribution of infected individuals among the different stages, and upper and
lower bounds of the final size of the epidemic. Some results on the total
number of infected, including the existence of an initial rising in disease prevalence, are
developed in Section \ref{sec4}. In Section \ref{sec5}, a probabilistic
approach to model the number of contacts allows the characteristics of this
number of contacts to be used explicitly in order to study their influence in
the model, which opens these models to different applications. We summarize
our results and perspectives in the final Discussion section.

\smallskip

\section{A staged progression epidemic model with general incidence}

\label{sec2}

This section introduces the formulation of a discrete-time staged progression
epidemic model with general incidence, no demography, and no re-infection, and
some of its basic properties are analyzed.

The staged progression (SP) epidemic model assumes
\cite{hyman_differential_1999,brauer_discrete_2010} a homogeneous
subpopulation of susceptibles, while the infected subpopulation is subdivided
into $n$ classes corresponding to different infection stages. When a
susceptible individual is infected, it enters the first of these infected
classes and progresses sequentially through the rest up to the $n$-th class,
from where it passes to the removed subpopulation, see Figure \ref{fig1}.

\setlength{\unitlength}{1.8ex} 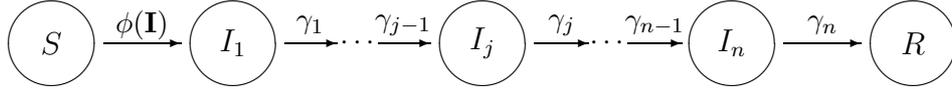
\begin{figure}[th]
\begin{center}
\begin{picture}(46,6)
\put(4,3){\circle{4}}
\put(12,3){\circle{4}}
\put(17.6,3){\makebox(0,0){$\cdots$}}
\put(23,3){\circle{4}}
\put(28.6,3){\makebox(0,0){$\cdots$}}
\put(34,3){\circle{4}}
\put(42,3){\circle{4}}
\put(4,3){\makebox(0,0){\large $S$}}
\put(12,3){\makebox(0,0){\large $I_1$}}
\put(23,3){\makebox(0,0){\large $I_j$}}
\put(34,3){\makebox(0,0){\large $I_n$}}
\put(42,3){\makebox(0,0){\large $R$}}
\put(6.3,3){\vector(1,0){3.4}}
\put(14.3,3){\vector(1,0){2.3}}
\put(18.4,3){\vector(1,0){2.3}}
\put(25.3,3){\vector(1,0){2.4}}
\put(29.4,3){\vector(1,0){2.4}}
\put(36.3,3){\vector(1,0){3.4}}
\put(8,3.8){\makebox(0,0){$\phi(\mathbf{I})$}}
\put(15.4,3.8){\makebox(0,0){$\gamma_1$}}
\put(19.5,3.8){\makebox(0,0){$\gamma_{j-1}$}}
\put(26.5,3.8){\makebox(0,0){$\gamma_{j}$}}
\put(30.6,3.8){\makebox(0,0){$\gamma_{n-1}$}}
\put(38,3.8){\makebox(0,0){$\gamma_n$}}
\end{picture}
\end{center}
\caption{Disease Flowchart of the staged progression epidemic model
(\ref{modn1}-\ref{modn4}).}%
\label{fig1}%
\end{figure}

Variables $S(t)$, $I_{j}(t)$, for $j=1,\ldots,n$, and $R(t)$ denote, as usual,
the size of the subpopulation of susceptibles, of the $j$-th infected class,
and of the removed subpopulation, respectively, at time $t=0,1,2,\dots$.

For each $j=1,\ldots,n$, individuals in class $I_{j}$ proceed to the next
class with probability $\gamma_{j}\in(0,1)$, and so they remain in this class
with probability $1-\gamma_{j}$.

We assume that the incidence over time step $(t,t+1]$, $S(t)-S(t+1)$, is given
by $\phi(\mathbf{I}(t))S(t),$ where $\phi$ is a certain function and
$\mathbf{I}(t)$ the population vector collecting the infected classes
\[
\mathbf{I}(t)=\left(  I_{1}(t),I_{2}(t),...,I_{n}(t)\right)  ^{T}
,\ t=0,1,2,\ldots
\]
Therefore, $\phi(\mathbf{I})$ represents the probability that a susceptible
individual is infected during the time step of the model and so $1-\phi
(\mathbf{I})$ is the probability of escaping infection.

The equations of the model are:
\begin{subequations}
\label{modn}%
\begin{align}
S(t+1)  &  =\left(  1-\phi(\mathbf{I}(t))\right)  S(t),\label{modn1}\\
I_{1}(t+1)  &  =\left(  1-\gamma_{1}\right)  I_{1}(t)+\phi(\mathbf{I}%
(t))S(t),\label{modn2}\\
I_{j}(t+1)  &  =\left(  1-\gamma_{j}\right)  I_{j}(t)+\gamma_{j-1}
I_{j-1}(t),\quad j=2,...,n,\label{modn3}\\
R(t+1)  &  =R(t)+\gamma_{n}I_{n}(t). \label{modn4}%
\end{align}

Let $\left\Vert \ast\right\Vert $ denote the 1-norm in the corresponding space
$\mathbf{R}^{m}$. The prevalence, i.e., the total infected population, that we
denote $Z(t)$, is
\end{subequations}
\[
Z(t):=I_{1}(t)+\cdots+I_{n}(t)=\left\Vert \mathbf{I}(t)\right\Vert ,\ t\geq0
\]

In model \eqref{modn} the total population remains constant, i.e.,
\[
S(t)+Z(t)+R(t)=N>0,\ t\geq0.
\]

\smallskip

Let us introduce some notation. If $\mathbf{x=}[x_{ij}]\mathbf{,y=}[y_{ij}%
]\in\mathbb{R}^{n\times m}$, we denote $\mathbf{x\geq y}$ (resp.
$\mathbf{x>y}$\textbf{) }to express that $x_{ij}\geq y_{ij}$ (resp.
$x_{ij}>y_{ij}$) for all $i=1,..,n,\ j=1,...,m$.

For each $S\in\lbrack0,N]$, we denote by $\mathcal{U}(S)$ the set of allowable
values for $\mathbf{I}$ when there are $S$ susceptibles, i.e.,
\begin{equation}
\mathcal{U}(S):=\left\{  \mathbf{I}\in\mathbb{R}_{+}^{n}:\ \left\Vert
\mathbf{I}\right\Vert +S\leq N\right\}  , \label{cat01}%
\end{equation}
and we define
\[
\mathcal{U}:=\mathcal{U}(0)=\left\{  \mathbf{I}\in\mathbb{R}_{+}%
^{n}:\left\Vert \mathbf{I}\right\Vert \leq N\right\}
\]

Next we are going to specify the hypotheses that the function $\phi$ must
satisfy so that $\phi(\mathbf{I})$ plays the role of force of infection, in
which the different infection stages act with their characteristic
infectivities. Each infection stage may or may not be infectious. Without loss
of generality we assume that the last infective stage is infectious, i.e., $R$
is the class that includes all post-infectious individuals.

Throughout this work we will assume function $\phi$ meets the following conditions:

\noindent\textbf{Hypothesis (H).}

\begin{enumerate}
\item[(i)] $\phi:\mathcal{U}\rightarrow\lbrack0,1)$, $\phi(\mathbf{0})=0$, and
$\phi\in C^{2}(\mathcal{U})$.

\item[(ii)] $\boldsymbol{\nabla}\phi(\mathbf{I})\geq\mathbf{0}$, for
$\mathbf{I}\in\mathcal{U}$, and $\mathbf{r}=\left(  r_{1},...,r_{n}\right)
:=\boldsymbol{\nabla} \phi(\mathbf{0})$ verifies $r_{n}>0$.

\item[(iii)] $\phi$ is concave-down in $\mathcal{U}$, i.e., $D^{2}%
\phi(\mathbf{I})$ is a negative semi-definite matrix, for $\mathbf{I}%
\in\mathcal{U}$.
\end{enumerate}

It is immediate to check that the previous hypotheses imply the following
inequalities
\begin{align}
\text{For all } \mathbf{I}  &  \in\mathcal{U},\ \phi(\mathbf{I})\leq
\mathbf{r\cdot I}.\label{cot1}\\
\text{For all } \mathbf{I}  &  \in\mathcal{U},\ \phi(\mathbf{I})\geq
\boldsymbol{\nabla} \phi(\mathbf{I})\mathbf{\cdot I}.\label{cot2}\\
\text{There exists }C  &  >0\text{ such that for all } \mathbf{I}%
\in\mathcal{U},\ \phi(\mathbf{I})\geq\mathbf{r} \cdot\mathbf{I}-C\left\Vert
\mathbf{I}\right\Vert _{2}^{2}, \label{cot3}%
\end{align}
where $\left\Vert \ast\right\Vert _{2} $ denotes the 2-norm.

\smallskip

Some examples of functions used in the literature in discrete epidemic models
that meet the previous hypothesis are
\begin{equation}
\phi(\mathbf{I})=1-\exp\left(  -\left(  \beta_{1}I_{1}+\cdots+\beta_{n}%
I_{n}\right)  \right)  ,\ \beta_{j}\geq0,\ j=1,...,n,\ \beta_{n}>0,
\label{finc1}%
\end{equation}
which is the standard choice in the literature
(\cite{brauer_discrete_2010,hernandez-ceron_discrete_2013,hernandez-ceron_reproduction_2013}%
);
\begin{equation}
\phi(\mathbf{I})=\beta_{1}I_{1}+\cdots+\beta_{n}I_{n},\ \beta_{j}%
\geq0,\ j=1,...,n,\ \beta_{n}>0,\ \beta_{1}+\cdots+\beta_{n}\leq1/N,
\label{finc2}%
\end{equation}
used (for $n=1$) in \cite{allen_discrete-time_1994}, and
\begin{align}
\phi(\mathbf{I})  &  =\theta_{1}\left(  1-\exp\left(  -\beta_{1}I_{1}\right)
\right)  +\cdots+\theta_{n}\left(  1-\exp\left(  -\beta_{n}I_{n}\right)
\right)  ,\ \label{finc3}\\
\beta_{j}  &  \geq0,\ j=1,...,n,\ \beta_{n}>0,\ \theta_{j}\in(0,1),\ \theta
_{1}+\cdots+\theta_{n}=1,\nonumber
\end{align}
used in \cite{van_den_driessche_disease_2019} in the context of a SEIRS model
$(n=2$) in which a fraction of the susceptible population have contacts with
$I$ individuals and the rest with $E$ individuals.

Note that $r_{j}$ characterizes the infectivity of class $I_{j}$ in the
neighbourhood of $0$, i.e, the infectivity \textquotedblleft up to first order
terms\textquotedblright. In particular, $r_{j}=0$ for a certain $j$ means
that, although class $I_{j}$ might be infectious, it is not infectious
\textquotedblleft up to first order terms\textquotedblright. Notices that in
the particular case of the functions $\phi$ given by (\ref{finc1}%
-\ref{finc3})), a class is infective if and only if is infective up to first
order terms.

By requiring $\mathbf{r}\geq\mathbf{0}$ and $r_{n}>0$, we are saying that the
infected classes may or may not be infectious except for the last last one,
$I_{n}$, that must necessarily be infectious up to first order. This does not
entail any loss of generality, for the last non-infectious classes can be
lumped into the removed class $R$.

Particular cases of these SP models are SIR ($n=1$), SEIR ($n=2$) and SE2I2R
models. In the latter, there are two compartments for each infective state,
i.e., the infective classes are $E_{1}$, $E_{2}$, $I_{1}$ and $I_{2},$ so that
the total length of time spent in classes $E$ and $I$ is the sum of two
geometric distributions.

\vspace{2ex}

%{\color{red} QUITAR y quizá citar en la introducción.
%A first example of model \eqref{modn} is a SEIR type model, in which
%$I_{1}=E$, $I_{2}=I$. Using an incidence function of the kind (\ref{finc1}),
%we have
%\[
%\phi(E,I)=1-\exp\left(  -\left(  \beta_{E}E+\beta_{I}I\right)  \right)
%,\ 0\leq\beta_{E},\ \beta_{I}>0.
%\]
%Another example is a SE2I2R model in which there are two compartments for each
%infective state, i.e., they are $E_{1}$, $E_{2}$, $I_{1}$ and $I_{2},$ so that
%the total length of time spent in classes $E$ and $I$ is the sum of two
%geometric distributions. Using an incidence function of the kind
%(\ref{finc1}), one has%
%\[
%\phi(E_{1},E_{2},I_{1},I_{2})=1-\exp\left(  -\left(  \beta_{E}\left(
%E_{1}+E_{2}\right)  +\beta_{I}\left(  I_{1}+I_{2}\right)  \right)  \right)
%,\ 0\leq\beta_{E},\ \beta_{I}>0
%\]
%}

\textit{Formulation of model \eqref{modn} as a dynamical system}.

\smallskip

As the equation for $R(t)$ is uncoupled from the rest, we consider the
dynamical system associated to the equations (\ref{modn1}, \ref{modn2},
\ref{modn3}). Let
\[
\mathcal{V}:=\left\{  (S,\mathbf{I}):\mathbf{I}\in\mathbb{R}_{+}^{n}%
,\ S\geq0,\ S+\left\Vert \mathbf{I}\right\Vert \leq N\right\}  ,
\]
and let us define the map $\mathbf{T}:\mathcal{V}\rightarrow\mathbb{R}^{n+1}$
as
\[
\mathbf{T}(S,\mathbf{I})=\left(  L(S,\mathbf{I}),\mathbf{G}(S,\mathbf{I}%
)\right)  =\left(  F(S,\mathbf{I}),G_{1}(S,\mathbf{I}),....,G_{n}%
(S,\mathbf{I})\right)  ,
\]
where
\begin{equation}%
\begin{array}
[c]{l}%
L(S,\mathbf{I})=\left(  1-\phi(\mathbf{I})\right)  S,\quad G_{1}%
(S,\mathbf{I})=\left(  1-\gamma_{1}\right)  I_{1}+S\phi(\mathbf{I}),\text{
and}\\
\rule{0ex}{3ex}G_{j}(S,\mathbf{I})=\left(  1-\gamma_{j}\right)  I_{j}%
+\gamma_{j-1}I_{j-1}\text{ for }j=2,\ldots,n.
\end{array}
\label{puf19}%
\end{equation}
Thus, equations (\ref{modn1}, \ref{modn2}, \ref{modn3}) can be expressed as
\[
\big(S(t+1),\mathbf{I}(t+1)\big)=\mathbf{T}\big(S(t),\mathbf{I}(t)\big),
\]
and also, separating the equation of the susceptible from those of the
infected, as
\begin{equation}
S\left(  t+1\right)  =L(S(t),\mathbf{I}(t)),\ \mathbf{I}\left(  t+1\right)
=\mathbf{G}(S(t),\mathbf{I}(t)),\ t\geq0. \label{puf21}%
\end{equation}

It is immediate to check that $\mathbf{T}(\mathcal{V})\subset\mathcal{V}$ and
so the dynamical system is well defined.

\smallskip

From now on we will assume that the initial population verifies
\begin{equation}
S(0)>0,\ \mathbf{I}(0)\in\mathcal{U}(S(0)),\ \mathbf{I}(0)\neq\mathbf{0},
\label{CI}%
\end{equation}
i.e., initially both the susceptible and the infected population are positive.

From \eqref{puf21} it is immediate to check that if (\ref{CI}) holds then
\begin{equation}
S(t)>0,\ \mathbf{I}(t)\in\mathcal{U}(S(0)),\ \mathbf{I}(t)\neq\mathbf{0}%
,\text{ for all }t\geq0, \label{modk25}%
\end{equation}
and,
\begin{equation}
\mathbf{I}(t)>\mathbf{0},\ R(t)>0,\text{ for all }t\geq n, \label{modk26}%
\end{equation}
that is, from $t=n$ onwards, every class has a positive size.

\vspace{2ex}

\textit{Long-term convergence of the variables of system \eqref{puf21}}

\smallskip

A simple consequence of \eqref{cot3} is that $\phi(\mathbf{I})\in(0,1)$ for
any $\mathbf{I}\in\mathcal{U}-\{\mathbf{0}\}$. Thus, for any initial condition
satisfying (\ref{CI}), sequence $S(t)$ is strictly decreasing and positive, so
there exists
\[
S_{\infty}:=\lim_{t\rightarrow\infty}S(t)\geq0.
\]
Adding equations (\ref{modn1}) and (\ref{modn2}) we obtain
\[
S(t+1)+I_{1}(t+1)= S(t)+I_{1}(t)-\gamma_{1}I_{1}(t),
\]
therefore, sequence $S(t)+I_{1}(t)$ is decreasing and positive, so convergent.
This implies that sequence $I_{1}(t)$ is also convergent and, taking limits in
the previous equality, we can conclude that
\[
\lim_{t\rightarrow\infty}I_{1}(t)=0.
\]
The previous argument can be repeated to reach the same conclusion with the
rest of the infection classes. By induction, for $j=2,3,\ldots$, let us assume
that $\lim_{t\rightarrow\infty}I_{i}(t)=0$ for $i<j$. Now, adding the $j+1$
first equations in \eqref{modn}, we have
\[
S(t+1)+I_{1}(t+1)+\cdots+I_{j}(t+1)= S(t)+I_{1}(t)+\cdots+I_{j}(t)-\gamma
_{j}I_{j}(t),
\]
that shows the convergence of sequence $S(t)+I_{1}(t)+\cdots+I_{j}(t)$ and so
that
\[
\lim_{t\rightarrow\infty}I_{j}(t)=0,\ j=2,...,k.
\]

We collect the results we have just proved in the following proposition.

\begin{proposition}
\label{prop:iinf} Every solution of system \eqref{puf21}, with initial
condition fulfilling (\ref{CI}), satisfies
\[
\lim_{t\rightarrow\infty}S(t)=S_{\infty},\ \lim_{t\rightarrow\infty}%
\mathbf{I}(t)=\mathbf{0},\ \lim_{t\rightarrow\infty}R(t)=N-S_{\infty
}.
\]

\end{proposition}

\smallskip

The following result guarantees that $S_{\infty}$ is positive.

\begin{proposition}
\label{prop:sinf} For every solution of system \eqref{puf21}, with initial
condition satisfying (\ref{CI}), we have $S_{\infty}>0$.
\end{proposition}

\begin{proof}
From (\ref{modn1}) we have
$$
S_{\infty}=S(0)\prod_{t=0}^{\infty}\Big(1-\phi\big(\mathbf{I}(t)\big) \Big)
$$
and so $S_{\infty}>0$ if and only if the infinite product is convergent, which
happens if and only if the infinite series of non-negative terms
$\sum^{\infty}_{t=0}\phi\big(\mathbf{I}(t)\big)$ converges. Using (\ref{cot1}) we have
\[
0\leq\overset{\infty}{\underset{t=0}{\sum}}\phi(\mathbf{I}(t))\leq\overset
{k}{\underset{j=1}{\sum}}r_{j}\overset{\infty}{\underset{t=0}{\sum}}I_{j}(t)
\]
and the series on the right converges because of Lemma \ref{res01} in the
Appendix. Therefore $S_{\infty}>0$ as we wanted to prove.
\end{proof}

\vspace{2ex}

\textit{Calculation of the basic reproduction number $\mathcal{R}_{0}$ of
system \eqref{puf21}}

\smallskip

We obtain $\mathcal{R}_{0}$ by applying the next-generation method as
developed in \cite{allen_basic_2008}. We naturally consider $I_{1},...,I_{n}$
as the infected states and $S$ as the only uninfected state.

%\begin{proposition}
%One has%
%\begin{equation}
%\mathcal{R}_{0}=N\delta\label{R0}%
%\end{equation}
%where
%\begin{equation}
%\delta:=\overset{k}{\underset{j=1}{\sum}}\frac{r_{j}}{\gamma_{j}}=\overset
%{k}{\underset{j=1}{\sum}}\frac{\frac{\partial\phi}{\partial I_{j}}%
%(\mathbf{0})}{\gamma_{j}}>0 \label{delta}%
%\end{equation}
%\end{proposition}
%\begin{proof}

Let $\mathbf{C}:=D_{\mathbf{I}}\mathbf{G}(N,\mathbf{0}),$ i.e., $\mathbf{C}$
is the matrix of the linearization of map $\mathbf{G}$ in \eqref{puf21} with
respect to the variables $\mathbf{I}$ evaluated at the disease free
equilibrium $(N,\mathbf{0})$.

After some calculations we obtain that $\mathbf{C}$ is equal to matrix
$\mathbf{B}(a)$ in (\ref{puf03}) (Lemma \ref{prop:lema10} in the Appendix) for
$a=N>0$, i.e., $\mathbf{C}=\mathbf{T}+\mathbf{F}(N),$ where matrix
$\mathbf{T}$ in the decomposition corresponds to transitions between infected
states, and $\mathbf{F}(N)$ corresponds to new infections. By definition,
$\mathcal{R}_{0}$ is the Net Reproductive Value \cite{li_applications_2002} of
matrix $\mathbf{C}$ corresponding to the aforementioned decomposition.
According to part (2) in Lemma \ref{prop:lema10} in the Appendix, we have%
\begin{equation}
\mathcal{R}_{0}=NRV(\mathbf{B}(N))=N\delta, \label{R0}%
\end{equation}
where
\begin{equation}
\delta:=\sum_{j=1}^{n}\frac{\partial\phi}{\partial I_{j}}(\mathbf{0}%
)/\gamma_{j}=\sum_{j=1}^{n}\frac{r_{j}}{\gamma_{j}}>0. \label{delta}%
\end{equation}

\smallskip

We can calculate $\mathcal{R}_{0}$ \eqref{R0} in the particular cases of
$\phi$ presented above. When $\phi$ is given by (\ref{finc1}) or (\ref{finc2})
one has $r_{j}=\beta_{j}$ and so
\[
\mathcal{R}_{0}=N\sum_{j=1}^{n}\frac{\beta_{j}}{\gamma_{j}}.
\]
This expression essentially matches (32) in \cite{brauer_discrete_2010}, where
$\phi$ is of type \eqref{finc1}. In the case of (\ref{finc3}), $r_{j}%
=\theta_{j}\beta_{j}$, and one obtains
\[
\mathcal{R}_{0}=N\sum_{j=1}^{n}\frac{\theta_{j}\beta_{j}}{\gamma_{j}}.
\]

\smallskip

\section{Asymptotic behavior of solutions}

\label{sec3}

{In this section we give some results regarding asymptotic properties of
solutions. On the one hand we study the asymptotic behavior of the components
of }$\mathbf{I}(t)$, characterizing their speed of convergence to zero and,
under some circumstances, their final {monotonicity. On the other hand, we
find some bounds for }$S_{\infty}$ which are related to the value of
$\mathcal{R}_{0}$ and show that, in the particular case of the incidence
function (\ref{finc1}), $\mathcal{R}_{0}$ is the solution of a certain
non-linear equations.

We first present a result that shows, loosely speaking, that all the
components of $\mathbf{I}(t)$ tend to zero at the same speed. In its statement
we work with a certain vector $\mathbf{v}$, which is defined from the
following matrix
\begin{equation}
\mathbf{A}:=\left(
\begin{array}
[c]{cccccc}%
1-\gamma_{1}+r_{1}S_{\infty} & r_{2}S_{\infty} & r_{3}S_{\infty} & \cdots &
r_{n-1}S_{\infty} & r_{n}S_{\infty}\\
\gamma_{1} & 1-\gamma_{2} & 0 & \cdots & 0 & 0\\
0 & \gamma_{2} & 1-\gamma_{3} & \cdots & 0 & 0\\
\vdots & \vdots & \vdots & \ddots & \vdots & \vdots\\
0 & 0 & 0 & \cdots & 1-\gamma_{n-1} & 0\\
0 & 0 & 0 & \cdots & \gamma_{n-1} & 1-\gamma_{n}%
\end{array}
\right)  . \label{modk23}%
\end{equation}
Since $S_{\infty}>0$ (Proposition \ref{prop:sinf}) and $r_{n}>0$ (Hypothesis
(H) (ii)), matrix $\mathbf{A}$ is irreducible and primitive, and so we can
define $\mathbf{v}$ to be the right Perron eigenvector of $\mathbf{A}$, i.e.,
the unique vector $\mathbf{v}=\left(  v_{1},...,v_{n}\right)  ^{T}$ that
verifies
\[
\label{punk10}\mathbf{Av}=\rho(\mathbf{A})\mathbf{v,\ v>0},\left\Vert
\mathbf{v}\right\Vert =1.
\]

\begin{proposition}
\label{prop31} Every solution of system \eqref{puf21}, with initial condition
fulfilling (\ref{CI}), satisfies that
\begin{equation}
\lim_{t\rightarrow\infty}\frac{\mathbf{I}(t)}{\left\Vert \mathbf{I}%
(t)\right\Vert }=\mathbf{v}. \label{modk41}%
\end{equation}
In particular, one has
\begin{equation}
\lim_{t\rightarrow\infty}\frac{I_{i}(t)}{I_{j}(t)}=\frac{v_{i}}{v_{j}%
},\ i,j=1,...,n. \label{modk42}%
\end{equation}

\end{proposition}

\begin{proof}
Let $(S(0),\mathbf{I}(0))$ satisfying (\ref{CI}) be fixed. We will now express
the dynamics of $\mathbf{I}(t)$ (\ref{puf21}) as a non-autonomous matrix
system. By applying the mean value theorem we have
\[
\phi(\mathbf{I})=\phi(\mathbf{0})+\nabla\phi(\mathbf{c(I)})\cdot\mathbf{I}=
\nabla\phi(\mathbf{c(I)})\cdot\mathbf{I,\ I}\in\mathcal{U},
\]
where $\mathbf{c}$ is a function such that for each $\mathbf{I}\in\mathcal{U}%
$, $\mathbf{c}(\mathbf{I})$ belongs to the segment that joins $\mathbf{0}$ and
$\mathbf{I}$ and, moreover, $\lim_{\mathbf{I}\rightarrow\mathbf{0}}%
\mathbf{c}(\mathbf{I})=\mathbf{0}$ and $\mathbf{c}(\mathbf{0})=\mathbf{0}$.
Therefore, if we define
\[
r_{j}(t)=\frac{\partial\phi}{\partial I_{j}}(\mathbf{c(I}(t)\mathbf{)}%
),\ j=1,...,n,
\]
and matrix
\begin{equation*}
\rule{-3ex}{0ex}\mathbf{A}(t):=\left(
\begin{array}
[c]{cccccc}%
1-\gamma_{1}+r_{1}(t)S(t) & r_{2}(t)S(t) & r_{3}(t)S(t) & \cdots &
r_{n-1}(t)S(t) & r_{n}(t)S(t)\\
\gamma_{1} & 1-\gamma_{2} & 0 & \cdots & 0 & 0\\
0 & \gamma_{2} & 1-\gamma_{3} & \cdots & 0 & 0\\
\vdots & \vdots & \vdots & \ddots & \vdots & \vdots\\
0 & 0 & 0 & \cdots & 1-\gamma_{n-1} & 0\\
0 & 0 & 0 & \cdots & \gamma_{n-1} & 1-\gamma_{n}%
\end{array}
\right)  , \label{modk22}%
\end{equation*}
then we can write the second equation of system  \ref{puf21} in the form
\begin{equation*}
\mathbf{I}(t+1)=\mathbf{A}(t)\mathbf{I}(t),
\end{equation*}
and so
\begin{equation*}
\mathbf{I}(t)=\Big(\prod_{s=0}^{t-1}\mathbf{A}(s)\Big)\mathbf{I}(0),
\end{equation*}
Note that, due to (\ref{modk25}) and Hypothesis (H) (ii), $S(t)>0$ and $r_j(t)\geq 0$ for $j=1,\ldots,n$ and $t\geq 0$, therefore, $\mathbf{A}(t)$ is nonnegative for $t\geq 0$. The fact that $\gamma_{j}\in(0,1),\ j=1,...,n$, implies that all entries in the main diagonal of matrices $\mathbf{A}(t)$, $t\geq 0$, are strictly positive, so we can conclude that they are allowable \cite{seneta_non-negative_1981}. Moreover, from Proposition \ref{prop:iinf} we have
$$
\lim_{t\rightarrow\infty}\mathbf{A}(t)=\mathbf{A},
$$
that is a primitive matrix. Now we use Theorem 3.6. from \cite{seneta_non-negative_1981} on strong ergodicity of non-negative matrix products. It implies that, when $t\to\infty$, every column in matrix
$\prod_{s=0}^{t}\mathbf{A}(s)$ tends to be proportional to $\mathbf{v}$, the right Perron eigenvector of matrix $\mathbf{A}$, and so (\ref{modk41}) follows. Expression (\ref{modk42}) readily follows from (\ref{modk41}).
\end{proof}

\smallskip

In the following theorem we obtain an upper bound of $S_{\infty}$ when
$\mathcal{R}_{0}>1$ and a lower bound if $\mathcal{R}_{0}<1$.

\begin{theorem}
\label{teorSinf} Let $S(t)$ be the susceptible population of a solution of
system \eqref{puf21} with initial condition fulfilling (\ref{CI}), and let
$S_{\infty}$ be its limit when $t\rightarrow\infty$. Then:

\begin{enumerate}
\item For $\mathcal{R}_{0}>1$,
\begin{equation}
\label{seir101}S_{\infty}\leq\frac{N}{\mathcal{R}_{0}}=\frac{1}{\delta},
\end{equation}
($\delta=\sum_{j=1}^{n}r_{j}/\gamma_{j}$ (\ref{delta})).

\item For $\mathcal{R}_{0}<1$, if the initial condition additionally verifies
$\mathbf{I}(0)=(I^{0},0,..,0)$ (i.e., all initial infected individuals are in
the first class), we have
\begin{equation}
\label{punk75}S_{\infty}\geq S(0)\frac{1-\delta(S(0)+I^{0})}{1-\delta S(0)}
=S(0)\frac{1-\mathcal{R}_{0}(S(0)+I^{0})/N}{1-\mathcal{R}_{0}S(0)/N}>0.
\end{equation}

\end{enumerate}
\end{theorem}

\begin{proof}
\begin{enumerate}
\item Summing the equations for $I_{j}(t),\ j=1,...,n$ (\ref{modn2},\ref{modn3})
one has
\[
Z\left(  t+1\right)  =Z(t)-\gamma_{n}I_{n}(t)+S(t)\phi(\mathbf{I}(t)).
\]
As $I_{n}(t)>0$ for $t\geq n$, we can write
\begin{equation}\label{za16}
Z\left(  t+1\right)  =Z(t)+D(t)I_{n}(t),
\end{equation}
with $D(t)=-\gamma_{n}+S(t)\phi(\mathbf{I}(t))/I_{n}(t),\ t\geq n$.
Using the mean value theorem, we have%
$$
\phi\left(  \mathbf{I}(t)\right)=\nabla\phi\left(
c\left(  \mathbf{I}(t)\right)  \right)  \cdot\mathbf{I}(t)
=\sum_{j=1}^{n}\frac{\partial\phi}{\partial I_{j}}\left(  c\left(
\mathbf{I}(t)\right)  \right)  I_{j}(t),
$$
where $c\left(  \mathbf{I}(t)\right)  $ is a point of the segment that joins
$\mathbf{0}$ and $\mathbf{I}(t)$. Since $\mathbf{I}(t)$ tends to $\mathbf{0}$
when $t\rightarrow\infty$, so does $c\left(\mathbf{I}(t)\right)$ and so we
have $\lim_{t\rightarrow\infty}\partial\phi/\partial I_{j}\left(
c\left(  \mathbf{I}(t)\right)  \right)  =\partial\phi/\partial I_{j}\left(  \mathbf{0}\right)  =r_{j},\ j=1,...,n$. Moreover, from Proposition
\ref{prop31} we have that $\lim_{t\rightarrow\infty}I_{j}(t)/I_{n}(t)=v_{j}/v_{n},\ j=1,...,n$, and so
$$
\lim_{t\rightarrow\infty}\phi\left(  \mathbf{I}(t)\right)/I_{n}(t)= \sum_{j=1}^{n}r_{j}v_{j}/v_{n}
$$
and, therefore,
\begin{equation}\label{za15}
\lim_{t\rightarrow\infty}D(t)=-\gamma_{n}+\frac{S_{\infty}}{v_{n}}\sum
_{j=1}^{n}r_{j}v_{j}.
\end{equation}
Now, note that matrix $\mathbf{A}$ in (\ref{modk23}) corresponds to matrix
$\mathbf{B}(a)$ in Lemma \ref{prop:lema10}  in the Appendix, for $a=S_{\infty}%
$, i.e., $\mathbf{A}=\mathbf{B}(S_{\infty})$. Since $S_{\infty}>0$ we can
apply the results thereof. In particular, from \eqref{ka03} and \eqref{ka04}, one has
\begin{equation}\label{za12}
\sign\big(\rho(\mathbf{A})-1\big)=\sign\big(S_{\infty}-\frac
{1}{\delta}\big),
\end{equation}
and
\begin{equation}\label{za13}
\sign\big(S_{\infty}-\frac
{1}{\delta}\big)=\sign\big(\frac{v_j}{v_n}-\frac{\gamma_n}{\gamma_j}\big),\ j=1,2,...,n-1.
\end{equation}
In order to prove that $S_{\infty}\leq 1/\delta$ we will proceed by contradiction.
Let us assume that $S_{\infty}>1/\delta$. Then, using (\ref{za15},
\ref{za12}, \ref{za13}) we have $v_{j}/v_{n}>\gamma_{n}/\gamma_{j}%
$,$\ j=1,...,n-1$ and so
$$
\lim_{t\rightarrow\infty}D(t)=-\gamma_{n}+S_{\infty}\sum_{j=1}^{n}r_{j}\frac{v_{j}}{v_{n}}
>-\gamma_{n}+\frac{1}{\delta}\sum_{j=1}^{n}r_{j}\frac{\gamma_{n}}{\gamma_{j}}=0,
$$
i.e., $\lim_{t\rightarrow\infty}D(t)>0$ and so there exists $t_{0}\geq0$ such
that $D(t)>0$ for $t\geq t_{0}$. Then from (\ref{za16}) it follows
$Z\left(  t+1\right)  >Z(t),\ t\geq t_{0}$, and so $Z(t)$ cannot converge to $0$ in contradiction to Proposition \ref{prop:iinf}.
\item We will use that if $a(t),$ $t=0,1,2,...$ are non-negative numbers, then
\[
\prod_{t=0}^{m}(1-a(t))\geq1-\overset{m}{\underset{t=0}{\sum}}a(t),\ m\geq0,
\]
Then, from (\ref{modn1}) we have%
\[
\frac{S_{\infty}}{S(0)}=\prod_{t=0}^{\infty}(1-\phi\left(  \mathbf{I}%
(t)\right)  )\geq1-\overset{\infty}{\underset{t=0}{\sum}}\phi\left(
\mathbf{I}(t)\right)  \geq1-\overset{n}{\underset{j=1}{\sum}}r_{j}%
\overset{\infty}{\underset{t=0}{\sum}}I_{j}(t)
\]
where we have used (\ref{cot1}). Now from (\ref{pum01})%
\[
\overset{\infty}{\underset{t=0}{\sum}}I_{j}(t)=\frac{1}{\gamma_{j}}\left(
S(0)-S_{\infty}+\overset{j}{\underset{i=1}{\sum}}I_{i}(0)\right)  =\frac
{1}{\gamma_{j}}\left(  S(0)-S_{\infty}+I^{0}\right)  ,\ j=1,...,n,
\]
and so
\[
\frac{S_{\infty}}{S(0)}\geq1-\overset{n}{\underset{j=1}{\sum}}\frac{r_{j}%
}{\gamma_{j}}\left(  S(0)-S_{\infty}+I^{0}\right)  =1-\delta\left(
S(0)-S_{\infty}+I^{0}\right)
\]
from where we have $S_{\infty}\left(  1-\delta S(0)\right) \geq S(0)(1-\delta\left(  S(0)+I^{0}\right))$. Using that $S(0)<N$ and that $\mathcal{R}_{0}<1$ we have $1-\delta S(0)>1-\delta N=1-\mathcal{R}_{0}>0$ and so \eqref{punk75}:
$$
S_{\infty}\geq S(0)\frac{1-\delta\left(  S(0)+I^{0}\right)}{1-\delta S(0)}
=S(0)\frac{1-\mathcal{R}_{0}\big(S(0)+I^{0}\big)/N}{1-\mathcal{R}_{0}S(0)/N}
$$
as desired. It is inmediate to check that the right hand side is a positive number.
\end{enumerate}
\end{proof}

\smallskip

In the case of incidence of exponential type \eqref{finc1} one can find an
equation whose solution is $S_{\infty}$. With its help, we prove that
inequality \eqref{seir101} is strict in this case. As we will see further on,
some results on the dynamics of the models are stronger whenever we can
guarantee that the inequality is strict.

\begin{proposition}
\label{eqSinf} Let us consider a solution of system \eqref{puf21} with initial
condition fulfilling (\ref{CI}) and incidence function (\ref{finc1}). Then,
$S_{\infty}$ is the only solution in $[0,S(0)]$ to the equation
\begin{equation}\label{modk30}
\log\frac{S(0)}{S_{\infty}}=\overset{n}{\underset{j=1}{\sum}%
}\frac{\beta_{j}}{\gamma_{j}}\Big(S(0)+\sum_{i=1}^{j}I_{i}(0)\Big)-\mathcal{R}%
_{0}\frac{S_{\infty}}{N}.
\end{equation}
Moreover,
\begin{equation}\label{modk31}
S_{\infty}<\frac{N}{\mathcal{R}_{0}}=\frac{1}{\delta}.
\end{equation}

\end{proposition}

\begin{proof}
We have,  using \eqref{modn1} and (\ref{pum01}),
\begin{multline*}
\frac{S_{\infty}}{S(0)}=\prod_{t=0}^{\infty}\big(1-\phi\left(\mathbf{I}(t)\right)\big) =\prod_{t=0}^{\infty}\exp\Big(-\overset{n}{\underset{j=1}{\sum}}\beta_{j}I_{j}(t)\Big)  \\=\exp\Big(-\overset{n}{\underset{j=1}{\sum}}\beta_{j}
\overset{\infty}{\underset{t=0}{\sum}}I_{j}(t)\Big)
=\exp\left(-\overset{n}{\underset{j=1}{\sum}}\frac{\beta_{j}}{\gamma_{j}}\Big(S(0)-S_\infty+
\sum_{i=1}^{j}I_i(0)\Big)\right)\\
=\exp\left(-\overset{n}{\underset{j=1}{\sum}}\frac{\beta
_{j}}{\gamma_{j}}\Big(S(0)+\sum_{i=1}^{j}I_i(0)\Big)\right)
\exp\left(\mathcal{R}_{0}\frac{S_{\infty}}{N}\right)
\end{multline*}
that leads to (\ref{modk30}).
\smallskip
To prove \eqref{modk31}, let us define function $g$ as
\begin{equation}\label{modk33}
g(x):=\log\frac{S(0)}{x}-\overset{n}{\underset{j=1}{\sum}}\frac{\beta
_{j}}{\gamma_{j}}\Big(S(0)+\sum_{i=1}^{j}I_i(0)\Big)+\mathcal{R}_{0}\frac{x}{N}.
\end{equation}
Clearly, $S_{\infty}$ is a solution to the equation $g(x)=0$ in $(0,S(0))$. From $g'(x)=-1/x+\delta$ we know that $g$ is decreasing on $(0,1/\delta)$ and increasing on $(1/\delta,S(0))$. Now, from $g(S(0))=-\sum_{j=1}^{n}\beta_{j}/\gamma_{j}\big(\sum_{i=1}^{j}I_i(0)\big)<0$ and $\lim_{x\to 0^+}g(x)=+\infty$, we deduce that $g(x)=0$ possesses just one solution in $(0,S(0))$, and this solution belongs to the interval $(0,1/\delta)$. As this solution must be $S_\infty$ we have that  $S_\infty<1/\delta$.
\end{proof}
If the initial condition verifies $\mathbf{I}(0)=(I^{0},0,..,0)$ and
$S(0)+I^{0}=N$, equation \eqref{modk30} simplifies to
\begin{equation}
\label{modk32}\log\frac{S(0)}{S_{\infty}}=\mathcal{R}_{0}\left(
1-\frac{S_{\infty}}{N}\right)  ,
\end{equation}
which coincides with equation (33) in \cite{brauer_discrete_2010}.

\smallskip

If we substitute the incidence function \eqref{finc1} for \eqref{finc2} in
Proposition \ref{eqSinf}, equation \eqref{modk30} is no longer satisfied but,
as the following result shows, inequality \eqref{modk31} is still true.

\begin{proposition}
\label{eqSinf1} Let us consider a solution of system \eqref{puf21} with
initial condition fulfilling (\ref{CI}) and let the incidence function be
(\ref{finc2}). Then,
\[
\label{modk31b}S_{\infty}<\frac{N}{\mathcal{R}_{0}}=\frac{1}{\delta}.
\]

\end{proposition}

\begin{proof}
We have,  using \eqref{modn1}, inequality
$1-x\leq e^{-x}$, and (\ref{pum01}):
\begin{multline*}
\frac{S_{\infty}}{S(0)}=\prod_{t=0}^{\infty}\Big(1-
\overset{n}{\underset{j=1}{\sum}}\beta_{j}I_{j}(t)\Big) \leq \prod_{t=0}^{\infty}\exp\Big(-\overset{n}{\underset{j=1}{\sum}}\beta_{j}I_{j}(t)\Big)  \\
=\exp\left(-\overset{n}{\underset{j=1}{\sum}}\frac{\beta
_{j}}{\gamma_{j}}\Big(S(0)+\sum_{i=1}^{j}I_i(0)\Big)\right)
\exp\left(\mathcal{R}_{0}\frac{S_{\infty}}{N}\right)
\end{multline*}
Using now the function $g$ \eqref{modk33} defined in the proof of Proposition \ref{eqSinf}, we have that $S_{\infty}$ satisfies $g(S_{\infty})>0$, and using the same reasoning as there we conclude that $S_{\infty}\in (0,1/\delta)$ as we wanted to prove.
\end{proof}

Our next proposition ensures that, whenever $S_{\infty}<1/\delta$, every
component $I_{j}(t)$ of $\mathbf{I}(t)$ is eventually strictly decreasing.
Before that, a result that ensures that, if for a specific time all the
infected classes reduce their number, the same happens from that moment onwards.

\begin{lemma}
\label{lem01} Let us consider system \eqref{puf21} with initial condition
fulfilling (\ref{CI}). If $t_{0}\geq0$ is such that $\mathbf{I}(t_{0}%
+1)<\mathbf{I}(t_{0})$, then $\mathbf{I}(t+1)<\mathbf{I}(t)$ for all $t\geq
t_{0}$.
\end{lemma}

\begin{proof} We have that $\mathbf{I}(t_0+1)<\mathbf{I}(t_0)$, i.e.,
$$
\mathbf{G}(S(t_0),\mathbf{I}(t_0))<\mathbf{I}(t_0).
$$
Checking the previous inequality for $(S(t_0+1),\mathbf{I}(t_0+1))$ we will have, by induction, proved the lemma. But this is immediate because $S(t_0+1)<S(t_0)$, $\phi(\mathbf{I}(t_0+1)\leq \phi(\mathbf{I}(t_0))$, $\mathbf{I}(t_0+1)<\mathbf{I}(t_0)$,  and $\mathbf{G}$ is strictly increasing in all its components:
$$
\mathbf{G}(S(t_0+1),\mathbf{I}(t_0+1))<\mathbf{G}(S(t_0),\mathbf{I}(t_0))=\mathbf{I}(t_0+1).
$$
\end{proof}

\begin{proposition}
\label{prop:monotonia} Let us consider a solution of system \eqref{puf21} with
initial condition fulfilling (\ref{CI}) and assume that $S_{\infty
}<N/\mathcal{R}_{0}=1/\delta$. Then, there exists $t_{0}\geq0$ such that
$\mathbf{I}(t+1)<\mathbf{I}(t)$ for all $t\geq t_{0}$.
\end{proposition}

\begin{proof}
With the help of Lemma \ref{lem01}, it is enough to find $t_0$ such that
$\mathbf{I}(t_0+1)<\mathbf{I}(t_0)$ to finish the proof.
$\mathbf{I}(t+1)<\mathbf{I}(t)$ is equivalent to the following set of inequalities
\begin{align*}
-\gamma_{1}I_{1}(t)+S(t)\phi(\mathbf{I}(t)) &  <0\\
-\gamma_{j}I_{j}(t)+\gamma_{j-1}I_{j-1}(t) &  <0,\ j=2,...,n.
\end{align*}
We will show that for large enough $t$ all of them can be met.
As $\phi(\mathbf{I})\leq \mathbf{r}\cdot \mathbf{I}$, to obtain $-\gamma_{1}I_{1}(t)+S(t)\phi(\mathbf{I})(t)<0$ we prove that
$-\gamma_{1}I_{1}(t)+S(t)\mathbf{r}\cdot \mathbf{I}(t)<0$, for $t$ large enough.
From the fact that $I_{n}(t)>0$ for all $t\geq n,$ we can divide $-\gamma
_{1}I_{1}(t)+S(t)\mathbf{r\cdot I}(t)$ by $I_{n}(t)$
$$
\frac{-\gamma_{1}I_{1}(t)+S(t)\mathbf{r\cdot I}(t)}{I_{n}(t)}=-\gamma_{1}%
\frac{I_{1}(t)}{I_{n}(t)}+S(t)\sum_{j=1}^{n}r_{j}\frac{I_{j}(t)}{I_{n}%
(t)},\ t\geq n,
$$
and passing to the limit $t\rightarrow\infty$ and using (\ref{modk42}) we
have
\begin{multline*}
L_1:= \lim_{t\rightarrow\infty}\frac{-\gamma_{1}I_{1}(t)+S(t)\mathbf{r\cdot I}(t)}{I_{n}(t)}
=-\gamma_{1}\frac{v_{1}}{v_{n}}+S_{\infty}\sum_{j=1}^{n}r_{j}\frac{v_{j}}{v_{n}}\\
=\frac{1}{v_{n}}\Big(  -\gamma_{1}v_{1}+S_{\infty}\sum_{j=1}^{n}r_{j}v_{j}\Big)
=\big(-\gamma_{1}v_{1}+S_{\infty}\mathbf{r}\cdot \mathbf{v}\big)/v_n.
\end{multline*}
Using Lemma \ref{prop:lema10} with $a=S_{\infty}$, i.e., $\mathbf{B}(S_{\infty})=\mathbf{A}$ \eqref{modk23}, and combining \eqref{ka03} and \eqref{ka04} we obtain
\begin{equation}\label{modk34}
\sign\big(S_{\infty}-\frac{1}{\delta}\big)
=\sign\big( -\gamma_{1}v_{1}+S_{\infty}\,\mathbf{r}\cdot\mathbf{v}\big)
=\sign\big(-\gamma_{j}v_{j}+\gamma_{j-1}v_{j-1}\big).
\end{equation}
This, together with assumption $S_{\infty}<1/\delta$, implies that limit $L_1$ is negative, and, thus, $-\gamma_{1}I_{1}(t)+S(t)\mathbf{r}\cdot \mathbf{I}(t)<0$ for $t$ large enough.
Now, let $j\in\left\{  2,...,n\right\}  $ be fixed. From the fact that
$I_{j-1}(t)>0$ for all $t\geq n$ we can divide $-\gamma_{j}I_{j}(t)+\gamma
_{j-1}I_{j-1}(t)$ by $I_{j-1}(t)$%
\[
\frac{-\gamma_{j}I_{j}(t)+\gamma_{j-1}I_{j-1}(t)}{I_{j-1}(t)}=-\gamma_{j}%
\frac{I_{j}(t)}{I_{j-1}(t)}+\gamma_{j-1},\ t\geq n,
\]
and passing to the limit $t\rightarrow\infty$ and using (\ref{modk42}) we
have
\[
L_{j}:=\lim_{t\rightarrow\infty}\frac{-\gamma_{j}I_{j}(t)+\gamma_{j-1}I_{j-1}
(t)}{I_{j-1}(t)}=\big(-\gamma_{j}v_{j}+\gamma_{j-1}v_{j-1}\big)/v_{j-1}.
\]
Again, we use \eqref{modk34} and $S_{\infty}<1/\delta$ to conclude that $L_{j}<0$.  This implies that $-\gamma_{j}I_{j}(t)+\gamma_{j-1}I_{j-1}(t)<0$
for large enough $t$, which concludes the proof.
\end{proof}

Note that Proposition \ref{eqSinf1} guarantees that for the classical incidence function \ref{finc1} one always has $S_{\infty}<1/\delta$ and from Proposition \ref{prop:monotonia} it follows that all the components of $\mathbf{I}(t)$ are eventually decreasing.

\smallskip

\section{Dynamics of the disease prevalence}

\label{sec4} In this section we study some properties connected with
the monotonicity of the prevalence of the disease, i.e., of the total infected
population $Z(t)$. In the first place we show that the behavior of $Z(t)$ can
be very different from the case of the SIR model, in which if the initial
infected population is small enough there are only two possible behaviors for
$Z(t)$. Moreover, we give some partial results regarding that behavior in the
general case and then stronger results whenever only
the last class $I_{n}$ is infectious.

In the classical SIR model, $\mathcal{R}_{0}<1$ implies that the prevalence decays monotonically to zero, whereas
$\mathcal{R}_{0}>1$ implies that for small enough values of the initial
infected population, the prevalence rises initially until it reaches a maximum
and from there it decays monotonically to zero.

In the general SP model \eqref{puf21}, Proposition \ref{prop:iinf} guarantees
that every infected class $I_{j}(t)$ tends to zero and, therefore, the same
happens to $Z(t)$.

Also, Proposition \ref{prop:monotonia} proves that if $S_{\infty}<1/\delta$,
from a certain time onwards all infected classes decay monotonically to zero.
When this occurs, obviously $Z(t)$ also decreases monotonically. Other than
that, the behaviour of $Z(t)$ is much more complex than in the SIR model.

\begin{figure}[h]
\centering
\begin{subfigure}[t]{0.45\textwidth}
\centering
\includegraphics[width=\linewidth]{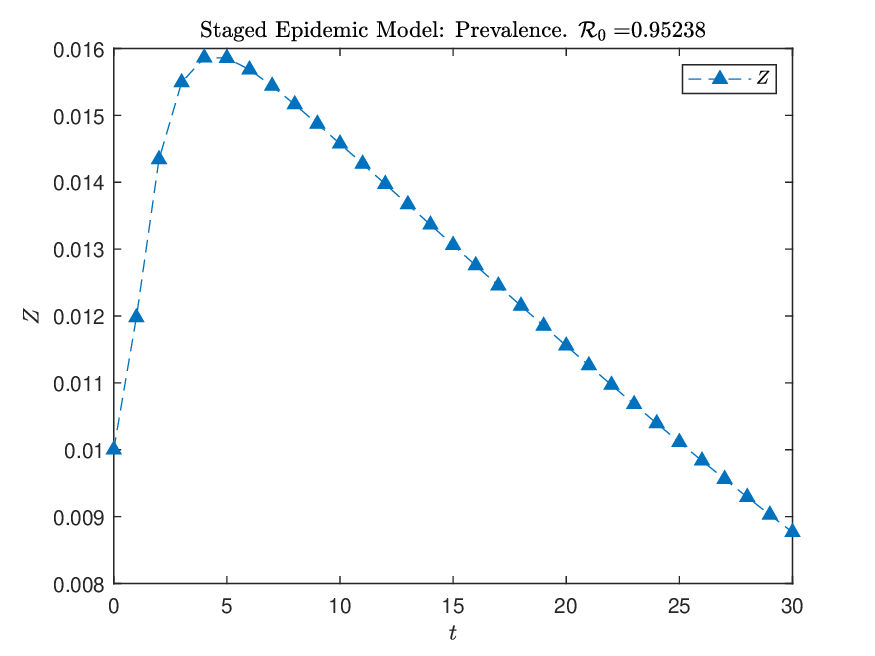}
\end{subfigure}
\qquad\begin{subfigure}[t]{0.45\textwidth}
\centering
\includegraphics[width=\linewidth]{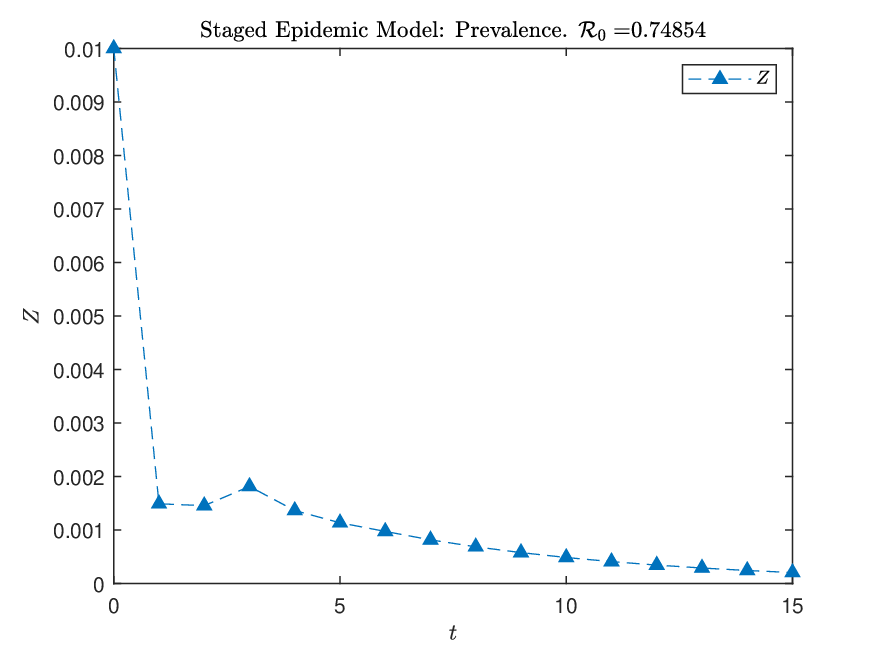}
\end{subfigure}
\caption{Two solutions of model \eqref{modn} for $n=3$ and the incidence function
(\ref{finc1}) in which $\mathcal{R}_{0}<1$ and however $Z(t)$ is not
monotonically decreasing. The parameter values are: \newline\textbf{Left}.
$N=1,\ \beta_{1}=\beta_{2}=0.2,\ \beta_{3}=0.1,\ \gamma_{1}=0.6,\ \gamma_{2}%
=0.7,\ \gamma_{3}=0.3,\ \mathbf{I}(0)=(0.01,0,0),\ R(0)=0$. \newline%
\textbf{Right}. $N=1,\ \beta_{1}=0.4,\ \beta_{2}=0.2,\ \beta_{3}=0.1,\ \gamma_{1}%
=\gamma_{3}=0.95,\ \gamma_{2}=0.9,\ \mathbf{I}(0)=(0,0,0.01),\ R(0)=0$.}%
\label{Fig1}%
\end{figure}

\begin{figure}[h]
\centering
\begin{subfigure}[t]{0.45\textwidth}
\centering
\includegraphics[width=\linewidth]{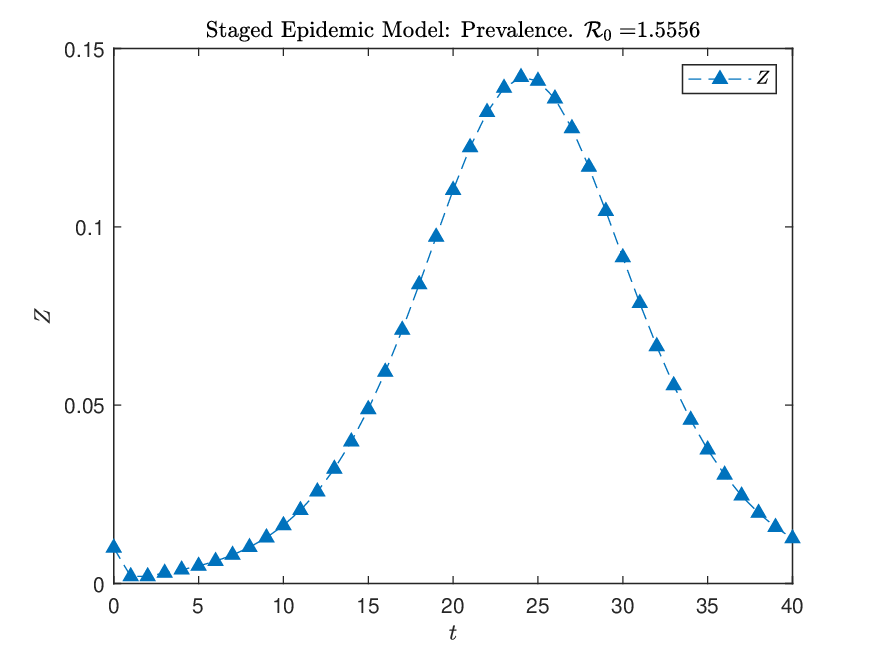}
\end{subfigure}
\qquad\begin{subfigure}[t]{0.45\textwidth}
\centering
\includegraphics[width=\linewidth]{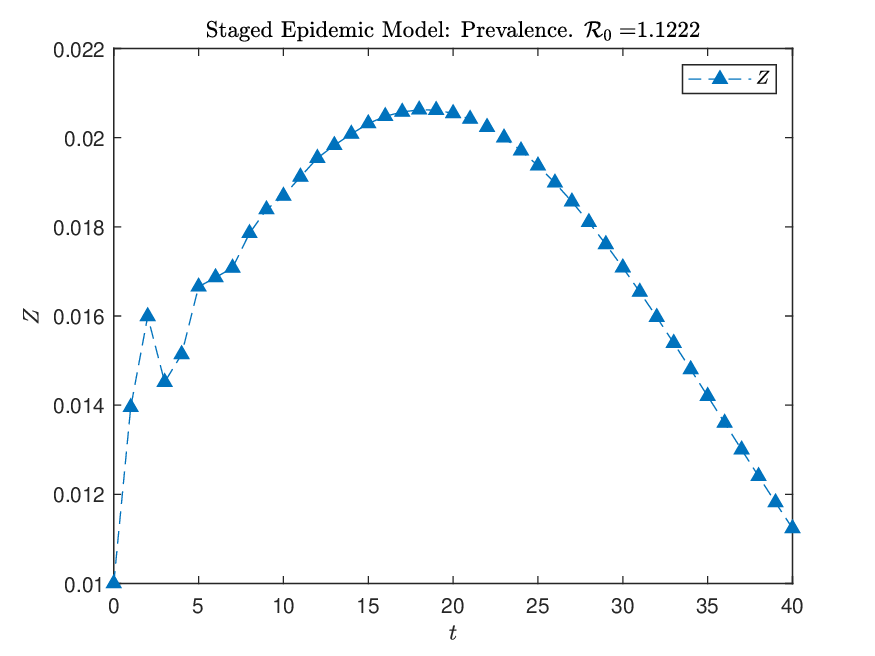}
\end{subfigure}
\par
\rule{0ex}{2ex}
\par
\centering
\begin{subfigure}[t]{0.45\textwidth}
\centering
\includegraphics[width=\linewidth]{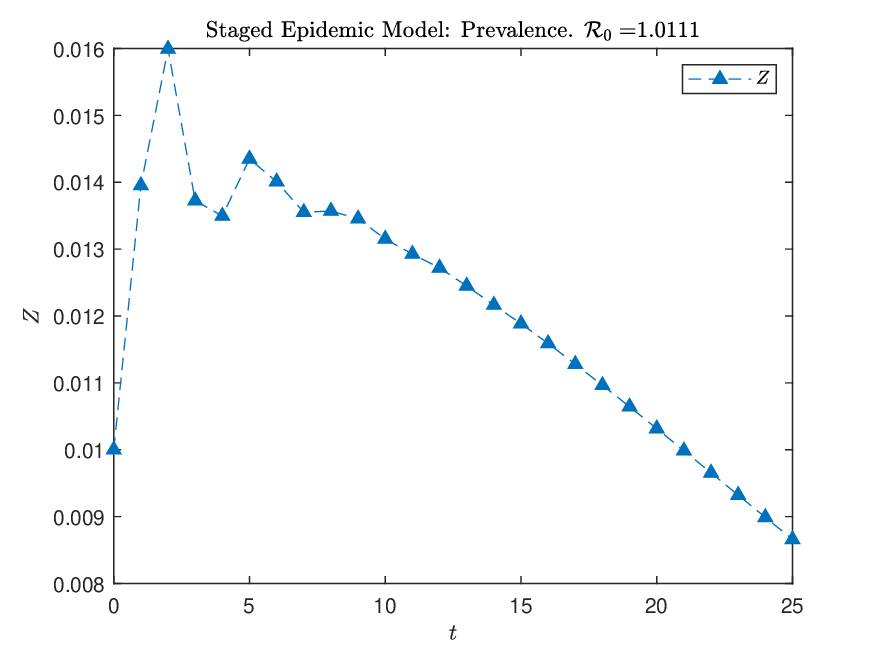}
\end{subfigure}
\caption{Three solutions of model \eqref{modn} for $n=3$ and the incidence function
(\ref{finc1}) in which $\mathcal{R}_{0}>1$ and however $Z(t)$ is not initially
increasing and then decreasing. The parameter values are: \newline\textbf{Top
Left}. $N=1,\ \beta_{1}=0.8,\ \beta_{2}=\beta_{3}=0.1,\ \gamma_{1}=0.6,\ \gamma_{2}%
=\gamma_{2}=0.9,\ \mathbf{I}(0)=(0,0,0.01),\ R(0)=0$. \newline\textbf{Top
Right}. $N=1,\ \beta_{1}=0.8,\ \beta_{2}=\beta_{3}=0.1,\ \gamma_{1}=0.6,\ \gamma
_{2}=\gamma_{2}=0.9,\ \mathbf{I}(0)=(0.01,0,0),\ R(0)=0$. \newline%
\textbf{Bottom}. $N=1,\ \beta_{1}=0.4,\ \beta_{2}=0.01,\ \beta_{3}=0.5,\ \gamma_{1}%
=\gamma_{2}=\gamma_{3}=0.9,\ \mathbf{I}(0)=(0,0,0.01),\ R(0)=0$. }%
\label{Fig2}%
\end{figure}

Figure \ref{Fig1} (resp. Figure \ref{Fig2}) show some solutions of model \eqref{modn} for $n=3$ and incidence function (\ref{finc1})in which $\mathcal{R}_{0}<1$ (resp. $\mathcal{R}_{0}>1$) and however the behavior of $Z(t)$ is not the one we described above regarding the SIR model.

The next propositions give some partial results on the monotonicity of $Z(t)$. In particular, we pay special attention to conditions under which the prevalence rises initially, i.e., $Z(1)>Z(0)$, and to those under which the prevalence decays monotonically to zero.
As it is the case in real epidemics, we are mainly interested in cases in
which the initial infected population is small with respect to $N$.

In order to simplify the notation, we define $\hat{\mathbf{I}}:=\left(
I_{1},...,I_{n-1}\right)  $, i.e., the vector with the first $n-1$ components
of $\mathbf{I}$. Therefore we have $\mathbf{I}=( \hat{\mathbf{I}},I_{n})$.
Recall that $\mathcal{U}(S)$, defined by (\ref{cat01}), is the set of
allowable $\mathbf{I}$ values when there are $S$ susceptibles ($\left\Vert
\mathbf{I}\right\Vert +S\leq N$).

We first consider the case in which at least one of the first $n-1$ infected
classes is infective up to the linear approximation.

\begin{proposition}
\label{prop:outbreak01} Let us consider system \eqref{puf21} with an initial
condition fulfilling (\ref{CI}). Let us assume that there exists $i\in\left\{
1,...,n-1\right\}  $ such that $r_{i}>0$ and, for some $t^{*}\geq0$,
$I_{i}(t^{*})>0$. Then there exists $c>0$ such that $0\leq I_{n}(t^{*})<c$
implies that $Z(t^{*}+1)>Z(t^{*})$.
\end{proposition}

\begin{proof}
Let $t\geq0$ be fixed. Summing equations (\ref{modn2}-\ref{modn3}) we have
\begin{equation}\label{za00}
Z(t+1)=Z(t)-\gamma_{n}I_{n}(t)+S(t)\phi\left(  \hat{\mathbf{I}}(t),I_{n}(t)\right).
\end{equation}
From the hipotheses on $\phi$, and the fact that $r_{i}=\partial\phi/\partial
I_{i}(\mathbf{0})>0$, we have that $\phi(\hat{\mathbf{I}}(t),I_{n}(t))>0$ whenever
$I_{i}(t)>0$. Therefore, using \eqref{za00}, we can assure that $Z(t^*+1)>Z(t^*)$ for every
$$
I_n(t^*)<c:=S(t^*)\phi\left(\hat{\mathbf{I}}(t^*),I_{n}(t^*)\right)/\gamma_{n}.
$$
\end{proof}

Note that, if there exists an infected class $I_{i},$ with $i\in\left\{
1,...,n-1\right\}  $, which is infective, up to the linear approximation,
$r_{i}>0$, and initially is not zero, $I_{i}(0)>0$, there are many
circumstances leading to the prevalence rising initially, $Z(1)>Z(0)$. It suffices to have no
infecteds in class $n$, i.e., $I_{n}(0)=0$, or an initial total infected
population sufficiently small, and it is independent of the initial value of
the susceptible population and of $\mathcal{R}_{0}$.

\smallskip

We treat now the case in which all the first $n-1$ infected classes are not
infective up to the linear approximation.

\begin{proposition}
\label{prop:outbreak02} Let us consider system \eqref{puf21} with an initial
condition fulfilling (\ref{CI}). Let us assume that $r_{1}=\cdots=r_{n-1}=0$.

\begin{enumerate}
\item If $t^{\ast}\geq0$ is such that $S(t^{\ast})<N/\mathcal{R}_{0}$ we have
$Z(t+1)<Z(t)$ for all $t\geq t^{\ast}$.

\item If $S(0)\in(N/\mathcal{R}_{0},N)$ (and this implies $\mathcal{R}_{0}>1$)
then there exists $\eta>0$ and a positive function $\alpha$ defined on
$(0,\eta)$ such that if $0<I_{n}(0)<\eta$ and $\big\Vert \mathbf{\hat{I}%
}(0)\big\Vert <\alpha(I_{n})$ then $Z(1)>Z(0)$.
\end{enumerate}
\end{proposition}

\begin{proof}
\begin{enumerate}
\item Since $r_{1}=\cdots=r_{n-1}=0$ we have $\mathcal{R}_{0}%
=Nr_{n}/\gamma_{n}$. Moreover, using (\ref{cot1}),
\begin{equation}\label{paf}
0\leq\phi(\mathbf{I}(t))\leq\mathbf{r}\cdot\mathbf{I}(t)=r_{n}I_{n}(t).
\end{equation}
Therefore, if $S(t^{\ast})<N/\mathcal{R}_{0}=\gamma_{n}/r_{n}$, one obtains
$$
S(t^{\ast})\phi(\mathbf{I}(t^{\ast}))-\gamma_{n}I_{n}(t^{\ast})<\frac
{\gamma_{n}}{r_{n}}r_{n}I_{n}(t^{\ast})-\gamma_{n}I_{n}(t^{\ast})=0,
$$
and so, according to (\ref{za00}), $Z(t^{\ast}+1)<Z(t^{\ast})$. As $S(t)$ is decreasing, $S(t)<N/\mathcal{R}_{0}$ for all $t\geq t^{\ast}$ and, reasoning as before,
it follows $Z(t+1)<Z(t)$ for all $t\geq t^{\ast}$.
\item If $I_n(0)>0$, then $\phi\left(\hat{\mathbf{I}}(0),I_{n}(0)\right)>0$, and, from \eqref{za00}, we can write
$$
Z(1)-Z(0)=\Big(S(0)-\frac{\gamma_nI_n(0)}{\phi\big(\hat{\mathbf{I}}(0),I_{n}(0)\big)}\Big)
\phi\big(\hat{\mathbf{I}}(0),I_{n}(0)\big),
$$
thus, the result is proved if we show that, for each $S(0)\in(N/\mathcal{R}_{0},N)$,
\begin{equation}\label{pach10}
\lim_{\mathbf{\hat{\mathbf{I}}}\rightarrow\mathbf{0}}\Big(\lim_{I_{n}\rightarrow 0^{+}}
\frac{\gamma_nI_n}{\phi\big(\hat{\mathbf{I}},I_{n}\big)}\Big)<S(0).
\end{equation}
Since $r_{n}=\partial\phi/\partial I_{n}(\mathbf{0})>0$ and $\partial\phi/\partial I_{n}$ is a continuous function, we have that $\partial\phi/\partial I_{n}(\mathbf{\hat{\mathbf{I}}},0)>0$ for $\big\Vert \mathbf{\hat{\mathbf{I}}}\big\Vert$ small enough. From \eqref{paf} we get that
$$
\lim_{I_n\to 0^+}\phi\big(\hat{\mathbf{I}},I_{n}\big)
=\phi\big(\hat{\mathbf{I}},0\big)=0,
$$
and, therefore,
$$
\lim_{\mathbf{\hat{\mathbf{I}}}\rightarrow\mathbf{0}}\Big(
\lim_{I_n\to 0^+}\frac{\gamma_nI_n}{\phi\big(\hat{\mathbf{I}},I_{n}\big)}\Big)=
\lim_{\mathbf{\hat{\mathbf{I}}}\rightarrow\mathbf{0}}\Big(
\frac{\gamma_n}{\partial\phi/\partial I_n\big(\hat{\mathbf{I}},0\big)}\Big)=\frac{\gamma_n}{r_n}=\frac{N}{\mathcal{R}_{0}}<S(0).
$$
\end{enumerate}
\end{proof}

As we see in item 1 of Proposition \ref{prop:outbreak02}, if $\mathcal{R}_{0}<1$ or
$S(0)\in(0,N/\mathcal{R}_{0})$, then we can take $t^{\ast}=0$ and so the prevalence tends to zero monotonically. If $S_{\infty}<N/\mathcal{R}_{0}$, which always occurs (Prop. \ref{eqSinf}) for the incidence functions \eqref{finc1} and \eqref{finc2}, then from Proposition
\ref{prop:monotonia} it follows that there exists $t^{*}$ from which the
total number of infecteds decreases monotonically towards zero. On the other
hand, item 2 shows that if $\mathcal{R}_{0}>1$ and $S(0)$ is large enough,
then there is an initial rising in the prevalence whenever the initial prevalence is
sufficiently small.

A particular case in which stronger results can be given is the situation in
which $I_{n}$ is the only infective class, i.e.,
\begin{equation}
\label{part01}\phi(\mathbf{I})=\varphi(I_{n})
\end{equation}
for a certain function $\varphi$. From Hypothesis (H) we have that
$\varphi:[0,N]\rightarrow\lbrack0,1)$ is a function such that\ $\varphi\in
C^{2}[0,N]$, $\varphi(0)=0$, $\varphi^{\prime}(x)\geq0$ and $\varphi
^{\prime\prime}(x)\leq0$ for all $x\in\lbrack0,N],$ and $r_{n}:=\varphi
^{\prime}(0)>0$. Note that $r_{1}=\cdots=r_{n-1}=0$ and so results in 2. of
Proposition \ref{prop:outbreak02} are applicable. For this particular case,
the incidence function of the form (\ref{finc2}) takes the form
\begin{equation}
\varphi(x)=\beta x,\ 0<\beta\leq1/N \label{finc4}%
\end{equation}
whereas (\ref{finc1}) and (\ref{finc3}) read
\begin{equation}
\varphi(x)=1-\exp\left(  -\beta x\right)  ,\ 0<\beta. \label{finc5}%
\end{equation}
In both cases one has $r_{n}=\beta$.

The simplest example of a model of this kind is the classical SEIR model in
which exposed individuals are not infectious.

\begin{proposition}
\label{prop:outbreak03} Let us consider system \eqref{puf21} with an initial
condition fulfilling (\ref{CI}). Let us assume that the only infective class
is $I_{n}$ and so function $\phi$ verifies (\ref{part01}). Then:

\begin{enumerate}
\item If $t^{\ast}\geq0$ is such that $S(t^{\ast})<N/\mathcal{R}_{0}$ we have
$Z(t+1)<Z(t)$ for all $t\geq t^{\ast}$.

\item If $S(0)\in(N/\mathcal{R}_{0},N)$ (and this implies $\mathcal{R}_{0}>1$)
then there exists $\eta>0$ such that if $0<I_{n}(0)<\eta$ then $Z(1)>Z(0)$.

\item If in addition to Hypothesis (H) funcion $\varphi$ verifies
\begin{equation}
\label{punk96}\varphi(x)\geq\frac{r_{n}x}{1+r_{n}x},\ x\in(0,N],
\end{equation}
then once sequence $Z(t)$ starts decreasing it decreases monotonically to
zero, i.e., if $t^{*}\geq0$ is such that $Z(t^{*}+1)<Z(t^{*})$, then
$Z(t+1)\leq Z(t)$ for all $t\geq t^{*}$.
\end{enumerate}
\end{proposition}

\begin{proof}
\begin{enumerate}
\item This is just a restatement of item 1. of Proposition
\ref{prop:outbreak02}.
\item This is essentially item 2. of Proposition \ref{prop:outbreak02}. We need to prove
\eqref{pach10}, that in this case simplifies to
$$
\lim_{\mathbf{\hat{\mathbf{I}}}\rightarrow\mathbf{0}}\Big(\lim_{I_{n}\rightarrow 0^{+}}
\frac{\gamma_nI_n}{\varphi(I_{n})}\Big)=
\lim_{I_{n}\rightarrow 0^{+}}
\frac{\gamma_nI_n}{\varphi(I_{n})}=\frac{\gamma_n}{r_n}=\frac{N}{\mathcal{R}_{0}}<S(0).
$$
\item Equation \eqref{za00} in this case reads as follows
\begin{equation}\label{za01}
Z(t+1)=Z(t)-\gamma_{n}I_{n}(t)+S(t)\varphi(I_{n}(t)).
\end{equation}
Let $t^*\geq0$ be such that $Z(t^*+1)<Z(t^*)$. From (\ref{za01}) it
follows that $I_{n}(t^*)>0$ (which implies that $\varphi(I_{n}(t))>0$ for
all $t\geq t^*$) and that
\begin{equation}\label{tan01}
S(t^*)<\gamma_{n}I_{n}(t^*)/\varphi(I_{n}(t^*)).
\end{equation}
Now, we want to prove that if (\ref{tan01}) holds, then $Z(t^*+2)<Z(t^*+1)$, i.e.,
\begin{equation}\label{tan02}
S(t^*+1)<\gamma_{n}I_{n}(t^*+1)/\varphi(I_{n}(t^*+1)).
\end{equation}
From (\ref{cot1}) we have $\varphi(I_n)\leq r_{n}I_n$ and so
\begin{equation}\label{tan03}
1/r_n\leq I_n/\varphi(I_n),
\end{equation}
and, from \eqref{punk96} we can obtain
\begin{equation}\label{tan04}
I_n/\varphi(I_n)\leq 1/r_n+I_n.
\end{equation}
Now, using \eqref{modn1}, \eqref{tan01}, \eqref{tan04}, and, finally \eqref{tan03}, we prove \eqref{tan02}:
\begin{multline*}
S(t^*+1)=(1-\varphi\big(I_{n}(t^*)\big))S(t^*)\\
<(1-\varphi\big(I_{n}(t^*)\big))\frac{\gamma_{n}I_{n}(t^*)}{\varphi\big(I_{n}(t^*)\big)}
=\gamma_{n}\Big(\frac{I_{n}(t^*)}{\varphi\big(I_{n}(t^*)\big)}-I_{n}(t^*)\Big)\\
<\gamma_{n}\Big(\frac{1}{r_n}+I_{n}(t^*)-I_{n}(t^*)\Big)=\frac{\gamma_{n}}{r_n}
\leq \frac{\gamma_{n}I_{n}(t^*+1)}{\varphi\big(I_{n}(t^*+1)\big)}.
\end{multline*}
Reasoning by induction, one has $Z(t+1)\leq Z(t)$
for all $t\geq t^*$.
\end{enumerate}
\end{proof}

In the case of incidence function (\ref{finc4}), condition (\ref{punk96})
holds trivially. In the case of (\ref{finc5}), expression (\ref{punk96}) is
equivalent to
\[
\beta x \leq\left(  1-\exp(-\beta x)\right)  \left(  1+\beta x\right)
,\ x\in(0,N],
\]
which, after some manipulations, reads
\[
1+\beta x \leq\exp(\beta x),\ \ x\in(0,N],
\]
which clearly holds using the property that $\exp(a)\geq1+a$ for all
$a\in\mathbf{R}$. Therefore, for these usual incidence functions whenever only the last class is infectious, once sequence $Z(t)$ starts decreasing it decreases monotonically to zero. Therefore, behaviors like those shown in Figure \ref{Fig2} are not possible.

\section{SP model with probability distribution of the number of contacts}
\label{sec5}

SP epidemic models have been used to introduce variability that basic
compartmental models lack. This variability is used, in general, to better
describe the dynamics of the disease, see \cite{hyman_differential_1999} with
a case of HIV transmission. In a second step, it also serves to analyze the
effect that treatments \cite{brauer_discrete_2010}, i.e., sanitary conditions,
or behavior changes \cite{brauer_simple_2011} may have on it. Obviously, we
can include in this second step the analysis of all the control measures such
as lockdowns or the use of masks, which have become frequent in times of pandemic.

Much of the behavioral changes or control measures related to the development
of the disease can be reflected within the model through contacts between
susceptible individuals and the pathogen.

In model \eqref{modn}, the function $\phi$ that models incidence does not
contemplate explicity the influence of contacts or, more specifically, of (a)
the contact probability distribution and (b) of the probability of infection
in each contact. However, in order to adequately reflect the changes in the
model induced by changes in the contacts, it is convenient that they appear
explicitly in the definition of the $\phi$ function. The purpose of this
section is to include in the general framework of model \eqref{modn}, the
effect of the points (a) and (b) above. In order to do so we will adapt what
is done in the discrete-time generic model presented in
\cite{seno_primer_2022}.

We first introduce the probability $p_{i}$ that a susceptible has $i$ contacts
($i=0,1,\ldots$) in a unit of time. Thus, we have the probability distribution
of the number of contacts $\{p_{i}\}$, with $\sum_{i=0}^{\infty}p_{i}=1$. We
now introduce the probability $\pi_{t}$ that a contact in the interval of time
$[t,t+1)$ between a susceptible and an infected individual produces an
infection. It seems clear that $\pi_{t}$ has to depend on the density of
infected in each of the stages at time $t$, i.e., of vector $\mathbf{I}(t)$.
Therefore, $\pi_{t}$ is a certain function $\Pi$ of $\mathbf{I}(t)$
\begin{equation}
\pi_{t}=\Pi\big(\mathbf{I}(t)\big). \label{s51}%
\end{equation}

The probability of a susceptible not being infected in interval $[t,t+1)$ if
there are $i$ contacts is $\big(1-\Pi\big(\mathbf{I}(t)\big)\big)^{i}$ and,
therefore, the probability of escaping infection in $[t,t+1)$ is
\[
\sum_{i=0}^{\infty}p_{i}\big(1-\Pi\big(\mathbf{I}(t)\big)\big)^{i}.
\]
To introduce this probability in model \eqref{modn}, it is enough to define
the function $\phi$ as
\begin{equation}
\phi\big(\mathbf{I}(t)\big)=1-\sum_{i=0}^{\infty}p_{i}\big(1-\Pi
\big(\mathbf{I}(t)\big)\big)^{i}. \label{s52}%
\end{equation}
The model is then as follows
\begin{subequations}
\label{modnp}%
\begin{align}
S(t+1)  &  =\Big(\sum_{i=0}^{\infty}p_{i}\big(1-\Pi\big(\mathbf{I}%
(t)\big)\big)^{i}\Big)S(t),\label{modnp1}\\
I_{1}(t+1)  &  =\left(  1-\gamma_{1}\right)  I_{1}(t)+\Big(1-\sum
_{i=0}^{\infty}p_{i}\big(1-\Pi\big(\mathbf{I}(t)\big)\big)^{i}%
\Big)S(t),\label{modnp2}\\
I_{j}(t+1)  &  =\left(  1-\gamma_{j}\right)  I_{j}(t)+\gamma_{j-1}%
I_{j-1}(t),\quad j=2,...,n,\label{modnp3}\\
R(t+1)  &  =R(t)+\gamma_{n}I_{n}(t). \label{modnp4}%
\end{align}

We want the function $\phi$ \eqref{s52} to satisfy Hypothesis (H), so that all
the results obtained for system \eqref{modn} in the previous sections can also
be applied to system \eqref{modnp}. To do this, it suffices to assume that
function $\Pi$ satisfies Hypothesis (H). Indeed, it is straightforward to
check that in this case also $\phi$ satisfies (H) for any probability
distribution of the number of contacts $\{p_{i}\}$.

It is clear that any system of the form \eqref{modnp} is a particular case of
system \eqref{modn} by choosing $\phi$ according to (\ref{s52}). Seen from the
other side, we can also consider any system \eqref{modn} as a particular case
of system \eqref{modnp}. Indeed, given a function $\phi$ it is enough to
define the function $\Pi$ of system \eqref{modnp} as $\Pi(\mathbf{I}%
):=\phi(\mathbf{I})$ and to choose as probability distribution of the number
of contacts $\{p_{i}\}$ the one corresponding to $p_{1}=1$ and $p_{i}=0$ for
all $i\neq1$.

If we use \eqref{R0} and \eqref{delta} to calculate the basic reproduction
number $\mathcal{R}_{0}$ of model \eqref{modnp} we get
\end{subequations}
\begin{equation}
\mathcal{R}_{0}=N\sum_{j=1}^{n}\frac{1}{\gamma_{j}}\,\frac{\partial\phi
}{\partial I_{j}}(\mathbf{0})=N\sum_{j=1}^{n}\frac{1}{\gamma_{j}}%
\Big(\sum_{i=1}^{\infty}i\,p_{i}\Big)\frac{\partial\Pi}{\partial I_{j}%
}(\mathbf{0})=N\bar{p}\sum_{j=1}^{n}\frac{1}{\gamma_{j}}\frac{\partial\Pi
}{\partial I_{j}}(\mathbf{0}), \label{s53}%
\end{equation}
where the total number of contacts is explicitly represented by its expected
value $\bar{p}$ according to probability distribution $\{p_{i}\}$.

A standard choice of the probability distribution $\{p_{i}\}$ is the Poisson
distribution, in which, for a fixed parameter $\lambda>0$ that represents its
mean value, the probability of $i$ contacts is $p_{i}=e^{-\lambda}\lambda
^{i}/i!$, and so
\begin{equation}
\phi\big(\mathbf{I}(t)\big)=1-\sum_{i=0}^{\infty}\frac{e^{-\lambda}\lambda
^{i}}{i!}\big(1-\Pi\big(\mathbf{I}(t)\big)\big)^{i}=1-\exp\left(  -\lambda
\Pi\big(\mathbf{I}(t)\big)\right)  . \label{s54}%
\end{equation}

In \eqref{s54} a large family of functions $\phi$ is defined, each one
corresponding to a choice of function $\Pi$. The mean number $\lambda$ of
contacts is specified in the expression. This should serve to study how
control measures or changes in behavior affect the development of a disease
that is representable in terms of contacts.

A possible generalization of the framework presented in this section is to
work with a probability distribution that distinguishes contacts according to
different stages of infection.

\section{Discussion}
\label{disc}

The importance of representing the variability of infectiousness in epidemic models is clear. SP models are a good tool to describe the variability due to the different states that an infected individual goes through. Continuous-time SP models have been used for this purpose for at least a couple of decades.

Discrete-time SP models, which were proposed in \cite{brauer_discrete_2010}, have hardly been used. This has been the case despite the aforementioned advantages of being easier to compare with the data and of a direct numerical  exploration.

In this work, a discrete-time SP model with a general incidence function has been analysed, from a theoretical point of view. The model does not include reinfection or demographic turnover, since its objective is not so much to anticipate long-term behavior as to provide information on the existence or not of epidemic outbreaks. In this sense, the results of the analysis of the model do not go through finding the attractors of the associated dynamical system, but rather focus on the study of issues such as the evolution over time of the disease prevalence, especially in its beginnings, and in the final size of the epidemic.

The consideration of a model of this type to represent the evolution of a disease has the clear purpose of serving as a test bench in which to contrast the different measures aimed at combating its progression.

Intervention policies in epidemic outbreaks include, in most cases, restrictions on contacts between individuals in one way or another. Sometimes this is caused by the individual's own behavioral change, and sometimes by decisions made by public health authorities. It is therefore interesting to be able to explicitly reflect the number of contacts in the model. This facilitates the evaluation of the impact of this type of intervention.
In Section \ref{sec5} an alternative formulation of the model has been proposed in which the probability distribution of the number of contacts appears explicitly.

The proposed model supports various natural extensions, which we have not addressed in order to obtain sufficiently clear results within a reasonable length. One of these extensions has to do with the possibility of amelioration. It consists of individuals not only advancing through the chain of infected states, but they can also go backwards. Another interesting extension is to contemplate, instead of one, two lines of advance of the infected individuals. These could represent, for infected individuals, those detected and those not detected or, also, those treated and those not treated.

Since SP models reflect the variability of infectiousness, another very important step is to treat them within the framework of meta-populations.
These meta-populations must be considered in a broad sense. The different items between which individuals transit could represent spatial places but also activities, behaviors, or any other classification that implies notable changes in infectiousness.
In these cases, it is not surprising that the processes, on the one hand, of progression through the infectious states and, on the other, of transition between the different items of the meta-population, can be considered acting on different time scales. We will address this issue in a future publication following the scheme developed in \cite{bravo_de_la_parra_reduction_2023}. In it, discrete-time models are proposed that combine two processes associated with different time scales, as well as the possible analysis of the model through a reduced model.

\section*{Funding}

Authors are supported by Ministerio de Ciencia e Innovaci\'on (Spain), Project PID2020-114814GB-I00.

\section{References}

\bibliographystyle{tfq}

\bibliography{MCMDS_DiscreteSPEpiMod_Sanz_Bravo}

\section{Appendix}

%\subsection{Some lemmas}

\begin{lemma}
\label{res01} Let us consider system \eqref{puf21} with an initial condition
fulfilling (\ref{CI}). For any $t_{0}\geq0$ we have%
\begin{equation}
\overset{\infty}{\underset{t=t_{0}}{\sum}}I_{j}(t)=\frac{1}{\gamma_{j}}\left(
S(t_{0})-S_{\infty}+\overset{j}{\underset{i=1}{\sum}}I_{i}(t_{0})\right)
,\ j=1,...,n, \label{pum01}%
\end{equation}
so that in particular
\[
\overset{\infty}{\underset{t=t_{0}}{\sum}}I_{n}(t)=\frac{1}{\gamma_{n}}\left(
S(t_{0})-S_{\infty}+\left\Vert I(t_{0})\right\Vert \right)  . \label{pum02}%
\]

\end{lemma}

\begin{proof}
From (\ref{modn1}-\ref{modn2}) we have
\[
I_{1}(\tau+1)-I_{1}(\tau)=-\gamma_{1}I_{1}(\tau)+\phi(\mathbf{I}(\tau
))S(\tau)=-\gamma_{1}I_{1}(\tau)+S(\tau)-S(\tau+1),\ \tau\geq0
\]
so that, summing from $\tau=t_{0}$ to $\tau=t-1,$
\[
I_{1}(t)-I_{1}(t_{0})=-\gamma_{1}\overset{t-1}{\underset{\tau=t_{0}}{\sum}%
}I_{1}(\tau)+S(t_{0})-S(t),\ t_{0}\geq0,\ t\geq1
\]
that, taking limits $t\rightarrow\infty$ and using that $\lim_{t\rightarrow
\infty}S(t)=S_{\infty}$ and $\lim_{t\rightarrow\infty}I_{j}(t)=0,\ j=1,...,n$, leads to
\begin{equation}
\gamma_{1}\overset{\infty}{\underset{t=t_{0}}{\sum}}I_{1}(t)=S(t_{0}%
)-S_{\infty}+I_{1}(t_{0}),\ t_{0}\geq0, \label{kan01}%
\end{equation}
which is (\ref{pum01}) for $j=1$. Now, from (\ref{modn3}),
\[
\gamma_{2}I_{2}(t)=\gamma_{1}I_{1}(t)+\left(  I_{2}(t)-I_{2}(t+1\right)  )
\]
so that, summing from $\tau=t_{0}$ to $\tau=t-1,$%
\[
\gamma_{2}\overset{t-1}{\underset{\tau=t_{0}}{\sum}}I_{2}(\tau)=\gamma
_{1}\overset{t-1}{\underset{\tau=t_{0}}{\sum}}I_{1}(\tau)+I_{2}(t_{0}%
)-I_{2}(t)
\]
and thus, taking limits $t\rightarrow\infty$ and using (\ref{kan01}),
\[
\gamma_{2}\overset{\infty}{\underset{t=t_{0}}{\sum}}I_{2}(t)=\gamma
_{1}\overset{\infty}{\underset{t=t_{0}}{\sum}}I_{1}(t)+I_{2}(t_{0}%
)=S(t_{0})-S_{\infty}+I_{1}(t_{0})+I_{2}(t_{0})
\]
which is (\ref{pum01}) for $j=2$. Proceeding analogously one obtains
(\ref{pum01}) for $j=3,...,n$.
\end{proof}

\smallskip

\begin{lemma}
\label{prop:lema10} Let $a>0$ and let $\mathbf{B}(a)$ be the non-negative
matrix defined by
\begin{equation}
\mathbf{B}(a):=\mathbf{T}+\mathbf{F}(a) \label{puf03}%
\end{equation}
where
\[
\mathbf{T:}=\left(
\begin{array}
[c]{cccccc}%
1-\gamma_{1} & 0 & 0 & \cdots & 0 & 0\\
\gamma_{1} & 1-\gamma_{2} & 0 & \cdots & 0 & 0\\
0 & \gamma_{2} & 1-\gamma_{3} & \cdots & 0 & 0\\
\vdots & \vdots & \vdots & \ddots & \vdots & \vdots\\
0 & 0 & 0 & \cdots & 1-\gamma_{n-1} & 0\\
0 & 0 & 0 & \cdots & \gamma_{n-1} & 1-\gamma_{n}%
\end{array}
\right)  \geq0,\ \mathbf{F}(a):=a\left(
\begin{array}
[c]{c}%
\mathbf{r}\\
\mathbf{0}\\
\vdots\\
\mathbf{0}%
\end{array}
\right)  \geq0
\]
Then:

\begin{enumerate}
\item Matrix $\mathbf{B}(a)$ is irreducible and primitive, and $\rho
(\mathbf{T})<1$.

\item The Net Reproductive Value of $\mathbf{B}(a)$ corresponding to the
decomposition (\ref{puf03}) is $NRV(\mathbf{B}(a))=a\delta$, where
$\delta=\sum_{j=1}^{n}r_{j}/\gamma_{j}$ (\ref{delta}). Moreover,
\begin{equation}
\label{ka03}\sign\big(\rho(\mathbf{B}(a))-1\big)=\sign\big(a-\frac{1}{\delta
}\big).
\end{equation}

\item Let $\mathbf{v}=\left(  v_{1},...,v_{n}\right)  ^{T}$ be the (only)
probability normed positive right eigenvector of $\mathbf{B}(a)$ associated to
$\rho(\mathbf{B}(a))$ (the right Perron eigenvector of $\mathbf{B}(a)$), and
let $i,j\in\left\{  1,...,n\right\}  $, with $i<j$. Then
\begin{equation}
\label{ka04}\sign\big(\rho(\mathbf{B}(a))-1\big)=\sign\big(\gamma_{i}%
v_{i}-\gamma_{j}v_{j}\big) =\sign\big( a\,\mathbf{r}\cdot\mathbf{v}-\gamma
_{1}v_{1}\big).
\end{equation}

\end{enumerate}
\end{lemma}

\begin{proof}
Let $a>0$ be fixed.
\begin{enumerate}
\item Since $\gamma_{i}$ and $1-\gamma_{i}$ belong to $(0,1)$ for $i=1,..,n$,
and $a,r_{n}>0$, the digraph associated to matrix $\mathbf{B}(a)$ is strongly
connected and so $\mathbf{B}(a)$ is irreducible. Moreover, the existence of
cycles of length 1 implies that $\mathbf{B}(a)$ is primitive.
As matrix $\textbf{T}$ is triangular, we have $\rho(\mathbf{T})=\max_{j=1,\ldots,n}\{1-\gamma_j\}<1$.
\item Once we have seen that $\rho(\mathbf{T})<1$, the Net Reproductive Value
\cite{li_applications_2002} of matrix $\mathbf{B}(a)$ corresponding to the decomposition
(\ref{puf03}) is $NRV(\mathbf{B}(a))=\rho\left(  \mathbf{Q}(a)\right)  $ where
$\mathbf{Q}(a)\mathbf{:}=\mathbf{F}(a)\left(  \mathbf{I}-\mathbf{T}\right)
^{-1}$. After some calculations we have%
$$\rule{-3ex}{0ex}
\mathbf{Q}=a\left(
\begin{array}
[c]{c}%
\mathbf{r}\\
0\\
\vdots\\
0
\end{array}
\right)  \left(
\begin{array}
[c]{cccc}
1/\gamma_{1} & 0 & \cdots & 0 \\
1/\gamma_{2} & 1/\gamma_{2} & \cdots & 0 \\
\vdots & \vdots & \ddots & \vdots \\
1/\gamma_{n} & 1/\gamma_{n} & \cdots & 1/\gamma_{n}
\end{array}
\right)  =a\left(
\begin{array}
[c]{cccc}
\overset{n}{\underset{j=1}{\sum}}r_{j}/\gamma_{j} & \overset
{n}{\underset{j=2}{\sum}}r_{j}/\gamma_{j} & \cdots & r_{n}/\gamma_{n}\\
0 & 0 & \cdots & 0\\
\vdots & \vdots & \ddots & \vdots\\
0 & 0 & \cdots & 0
\end{array}
\right)
$$
and so $NRV(\mathbf{B}(a))=a\delta$.
We know that $\sign\big(\rho(\mathbf{B}(a))-1\big)=\sign\big(NRV(\mathbf{B}(a))-1\big)$
\cite{li_applications_2002}, so it follows \eqref{ka03}.
\item Let $\lambda:=\rho(\mathbf{B}(a))$. If we write $\mathbf{B}(a)\mathbf{v=}%
\lambda\mathbf{v}$ componentwise, we have
\begin{align*}
\left(  1-\gamma_{1}+ar_{1}\right)  v_{1}+ar_{2}v_{2}+\cdots+ar_{n}v_{n} &
=\lambda v_{1},\\
\gamma_{j-1}v_{j-1}+\left(  1-\gamma_{j}\right)  v_{j} &  =\lambda
v_{j},\ j=2,...,n,
\end{align*}
that can be expressed as
\begin{align*}
-\gamma_{1}v_{1}+ a\,\mathbf{r}\cdot\mathbf{v} &
=(\lambda-1) v_{1},\\
\gamma_{j-1}v_{j-1}-\gamma_{j}v_{j}=(\lambda-1)v_{j} &  =(\lambda-1)
v_{j},\ j=2,...,n.
\end{align*}
From the first equality it follows the equality between the first and the third term in \eqref{ka04}. To obtain the equality between the two first terms in \eqref{ka04}, it is enough to write the sum of the second equalities from $i$ to $j-1$:
$$
\gamma_{i}v_{i}-\gamma_{j}v_{j}=\sum_{k=i}^{j-1}\big(\gamma_{k}v_{k}-\gamma_{k+1}v_{k+1}\big)
=(\lambda-1)\sum_{k=i}^{j-1}v_{k}.
$$
\end{enumerate}
\end{proof}

\end{document}